\documentclass[10pt]{article}

\usepackage[a4paper]{geometry}
\usepackage{indentfirst}
\usepackage{hyperref}
\usepackage{amssymb, amsmath, amsthm}
\usepackage{multirow}
\usepackage[labelsep=none]{caption}
\usepackage[labelformat=simple]{subcaption}

\hypersetup{bookmarksnumbered, pdfstartview=FitH, pdfpagemode=UseNone}

\newcommand{\bfa}{\mathbf{a}}

\newcommand{\bfe}{\mathbf{e}}
\newcommand{\bff}{\mathbf{f}}
\newcommand{\bfx}{\mathbf{x}}

\newcommand{\bfz}{\mathbf{z}}
\newcommand{\sfM}{\mathsf{M}}
\newcommand{\sfO}{\mathsf{O}}
\newcommand{\sfZ}{\mathsf{Z}}
\newcommand{\shift}{{\lozenge}}
\DeclareMathOperator{\cf}{CF}

\newcommand{\mynewtheorem}[2]{\newtheorem{#1}{\indent #2}}
\newcommand{\myalttheorem}[2]{\newtheorem*{#1}{\indent #2}}
\mynewtheorem{lemma}{Lemma}
\mynewtheorem{proposition}{Proposition}
\myalttheorem{proposition*}{Proposition}
\mynewtheorem{theorem}{Theorem}
\myalttheorem{theorem*}{Theorem}
\newenvironment{myproof}[1][Proof]{\begin{proof}[\indent #1]}{\end{proof}}

\begin{document}

\title{\textbf{Powers of 2 in High-Dimensional Lattice Walks}}
\author{Nikolai Beluhov}
\date{}

\maketitle

\begin{abstract} Let $W_d(n)$ be the number of $2n$-step walks in $\mathbb{Z}^d$ which begin and end at the origin. We study the exponent of $2$ in the prime factorisation of this number; i.e., $w_d(n) = \nu_2(W_d(n))$. We show that, for each $d$, there is a relationship between $w_d(n)$ and the number $s_2(n)$ of $1$s in the binary expansion of $n$. For example, $w_d(n) = s_2(n)$ if $d$ is odd and $w_d(n) = 2s_2(n)$ if $\nu_2(d) = 1$; while $w_d(n) \ge 3s_2(n)$ if $\nu_2(d) = 2$. The pattern changes further when $\nu_2(d) \ge 3$. However, for each $d$, we give the best analogous estimate of $w_d(n)$ together with a description of all $n$ where equality is attained. The methods we develop apply to a wider range of problems as well, and so might be of independent interest. \end{abstract}

\section{Introduction} \label{intro}

Fix a positive integer $d$. For each positive integer $n$, let $W_d(n)$ be the number of $2n$-step walks in the $d$-dimensional integer lattice $\mathbb{Z}^d$ which begin and end at the origin. We assume that each step is of unit length; so, from each point of $\mathbb{Z}^d$, there are $2d$ options for our next move.

Then $W_1(n) = \binom{2n}{n}$ and $W_2(n) = \binom{2n}{n}^2$. When $d \ge 3$, though, closed-form formulas are no longer readily apparent. The most natural expression we can write down is perhaps as follows: For each basis vector $\bfe_i$, we must make the same number of steps pointing in direction $\bfe_i$ as in direction $-\bfe_i$. Denote this number by $x_i$. Then \[W_d(n) = \sum_{x_1 + x_2 + \cdots + x_d = n} \binom{2n}{x_1, x_1, x_2, x_2, \ldots, x_d, x_d}. \tag{\textbf{W}}\]

Or, alternatively, let $S_d$ be the Laurent polynomial $z_1 + 1/z_1 + z_2 + 1/z_2 + \cdots + z_d + 1/z_d$. Then $W_d(n)$ equals the constant term of $S_d^{2n}$; so, for each fixed $d$, we get what is known as a constant-term sequence. Other formulas as well as numerical data can be found on the On-Line Encyclopedia of Integer Sequences -- henceforth, OEIS. The entry \cite{O1} collects the full two-parameter array where both of $d$ and $n$ vary; and it links back to the individual entries for the sequences obtained with $d = 1$, $2$, $\ldots$, $6$.

Consider now the subset of those walks where the first $n$ steps are all positive and the latter $n$ steps are all negative. (We call a step ``positive'' when its direction is among $\bfe_1$, $\bfe_2$, $\ldots$, $\bfe_d$ and ``negative'' when it is among $-\bfe_1$, $-\bfe_2$, $\ldots$, $-\bfe_d$.) Let the number of such walks be $W^\star_d(n)$. Clearly, $W_d(n) = \binom{2n}{n}W^\star_d(n)$ and \[W^\star_d(n) = \sum_{x_1 + x_2 + \cdots + x_d = n} \binom{n}{x_1, x_2, \ldots, x_d}^2. \tag{$\mathbf{W}^\star$}\]

The numbers $W^\star_d(n)$ appear also in the combinatorics of words. An \emph{Abelian square} is a word of even length whose first and second half are anagrams \cite{RS}; so $W^\star_d(n)$ is the number of Abelian squares of length $2n$ over an alphabet with $d$ letters. The full two-parameter array is collected in the OEIS entry \cite{O2}, alongside links to the individual entries for the sequences obtained with $d = 1$, $2$, $\ldots$, $6$. The case of $d = 4$ is especially significant \cite{O3} -- the terms of the associated sequence are known as the \emph{Domb numbers}.

It was conjectured by Bhattacharya that the number $B(m, n)$ of $2m \times 2n$ balanced grid colourings exhibits certain arithmetical properties; see \cite{B}. A simple bijective argument shows that $B(1, n) = W_1(n)$ and $B(2, n) = W_3(n)$. (This seems to be a ``happy coincidence'', in the sense that $B(m, n)$ and $W_d(n)$ are not related in an obvious manner otherwise.) Hence, these arithmetical properties ought to be shared by $W_1(n)$ and $W_3(n)$ as well. Specifically, the implication is that the exponent of $2$ in the prime factorisation of each one of $W_1(n)$ and $W_3(n)$ must equal the number $s_2(n)$ of $1$s in the binary expansion of $n$.

It is natural now to wonder whether the numbers $W_d(n)$ exhibit similar arithmetical properties for any other values of $d$. A bit of thought quickly reveals that, in fact, the exponent of $2$ in $W_d(n)$ equals $s_2(n)$ for all odd $d$. By thinking a little harder, we can convince ourselves also that the exponent of $2$ in $W_d(n)$ equals $2s_2(n)$ for all $d$ with $d \equiv 2 \pmod 4$. We spell out both arguments in Section \ref{init}.

This feels like the beginning of a pattern. Given a positive integer $x$, let $\nu_2(x)$ denote the exponent of $2$ in the prime factorisation of $x$, and set $w_d(n) = \nu_2(W_d(n))$. The smallest value of $d$ we have not resolved yet is $d = 4$. The numerical data available on OEIS allows us to plot $w_4(n)$ over $n = 1$, $2$, $\ldots$, $50$. (Table \ref{w4t}.) However, at this point we suffer a temporary setback. Our naive expectation, extrapolating based on the families $d \equiv 1 \pmod 2$ and $d \equiv 2 \pmod 4$, would have been for $w_d(n)$ to equal $3s_2(n)$ if $d \equiv 4 \pmod 8$. This is, alas, not what we actually observe when $d = 4$.

\begin{table}[ht] \centering \footnotesize

\begin{tabular}{c|c|c|c|c|c|c|c|c|c|c|c|c|c|c|c|c|}
\hline
$n$ & 1 & 2 & 3 & 4 & 5 & 6 & 7 & 8 & 9 & 10 & 11 & 12 & 13 & 14 & 15 & 16\\
\hline
$s_2(n)$ & 1 & 1 & 2 & 1 & 2 & 2 & 3 & 1 & 2 & 2 & 3 & 2 & 3 & 3 & 4 & 1\\
\hline
$w_4(n)$ & 3 & 3 & 10 & 3 & 6 & 8 & 12 & 3 & 6 & 6 & 11 & 8 & 11 & 12 & 16 & 3\\
\hline
\end{tabular}

\bigbreak

\begin{tabular}{c|c|c|c|c|c|c|c|c|c|c|c|c|c|c|c|c|}
\hline
$n$ & 17 & 18 & 19 & 20 & 21 & 22 & 23 & 24 & 25 & 26 & 27 & 28 & 29 & 30 & 31 & 32\\
\hline
$s_2(n)$ & 2 & 2 & 3 & 2 & 3 & 3 & 4 & 2 & 3 & 3 & 4 & 3 & 4 & 4 & 5 & 1\\
\hline
$w_4(n)$ & 6 & 6 & 12 & 6 & 9 & 13 & 17 & 8 & 11 & 11 & 16 & 12 & 15 & 16 & 20 & 3\\
\hline
\end{tabular}

\bigbreak

\begin{tabular}{c|c|c|c|c|c|c|c|c|c|c|c|c|c|c|c|c|c|c|}
\hline
$n$ & 33 & 34 & 35 & 36 & 37 & 38 & 39 & 40 & 41 & 42 & 43 & 44 & 45 & 46 & 47 & 48 & 49 & 50\\
\hline
$s_2(n)$ & 2 & 2 & 3 & 2 & 3 & 3 & 4 & 2 & 3 & 3 & 4 & 3 & 4 & 4 & 5 & 2 & 3 & 3\\
\hline
$w_4(n)$ & 6 & 6 & 14 & 6 & 9 & 11 & 15 & 6 & 9 & 9 & 14 & 13 & 16 & 17 & 21 & 8 & 11 & 11\\
\hline
\end{tabular}

\caption{} \label{w4t} \end{table}

Looking at our experimental data somewhat more closely, we find that $w_4(n)$ and $3s_2(n)$ do seem to be related after all -- though by means of an inequality instead of an identity, with $w_4(n) \ge 3s_2(n)$. Furthermore, the exact relation $w_4(n) = 3s_2(n)$ does hold for a lot of the values of $n$ in Table \ref{w4t}, even if it is not universal. By feeding these values back into the OEIS search engine, we find that they coincide with the first few positive integers whose binary expansions do not contain two adjacent $1$s. We will show, in Section \ref{proof-i}, that both of these observations generalise to all $d \equiv 4 \pmod 8$ and $n$.

The family $d \equiv 8 \pmod{16}$ comes next. Once again, there is an obvious extrapolation based on our current knowledge: $w_d(n) \ge 4s_2(n)$. Once again, it is thwarted by the experimental data. When we go past $\nu_2(d) = 2$, the slope of the associated linear function of $s_2(n)$ stabilises and remains equal to $3$ for all $\nu_2(d) \ge 3$; but we get a constant term, too, which does grow together with $\nu_2(d)$. The sequence of those $n$ for which our linear lower bound is exactly equal to $w_d(n)$ stabilises as well; however, it settles into a different pattern from the one which occurs when $\nu_2(d) = 2$. Our findings can be summarised as follows:

\begin{theorem*} The exponent $w_d(n)$ of $2$ in the prime factorisation of the number of $d$-dimensional $2n$-step origin-to-origin lattice walks $W_d(n)$ satisfies the following identities and inequalities:

(a) If $\nu_2(d) = 0$, then $w_d(n) = s_2(n)$.

(b) If $\nu_2(d) = 1$, then $w_d(n) = 2s_2(n)$.

(c) If $\nu_2(d) = 2$, then $w_d(n) \ge 3s_2(n)$. Equality is attained if and only if the binary expansion of $n$ does not contain two adjacent $1$s.

(d) If $\nu_2(d) \ge 3$, let $\delta = \nu_2(d) - 2$. Then $w_d(n) \ge 3s_2(n) + \delta$. Equality is attained if and only if, in the binary expansion of $n$, all of the $1$s come before all of the $0$s. \end{theorem*}

The linear lower bounds in parts (c) and (d) are the best possible, in the sense that neither the slope can be improved; nor, given the optimal slope, can the constant term be improved. This follows by the descriptions of the equality cases, and it remains true even if we restrict consideration to those $n$ for which both of $n$ and $s_2(n)$ are greater than some fixed quantity.

The equality conditions of parts (c) and (d) can also be rephrased in terms of forbidden subwords. Specifically, in part (c) equality is attained if and only if the binary expansion of $n$ does not contain the subword $11$; while in part (d) the forbidden subword becomes $01$. This hints at a connection to finite automata, which we make explicit in Sections \ref{abac-i}--\ref{abac-ii}.

We give two different proofs of the Theorem. For the first one, we develop a general technique which applies also to a range of other problems where we want some sum to be divisible by some prime power. Roughly speaking, the key idea is to split the sum into sub-sums so that each sub-sum behaves less like a sum and more like a product. This technique makes it possible to treat all values of $d$ in parts (c) and (d) ``uniformly'', and the resulting treatment is fully self-contained.

By contrast, our second proof is not self-contained; instead, it relies on one identity between hypergeometric sums due to Chan and Zudilin \cite{CZ}. We use this identity to resolve the case of $d = 4$. The remaining values of $d$ in parts (c) and (d) are then handled in several ``batches'' by means of a series of varied arguments. The whole series is ultimately rooted at the case when $d = 4$, indicating that most of the difficulty is concentrated there.

One application of the first proof's technique is particularly amusing, and we give it here by way of an example. For convenience, let $s = s_2(n)$. It is well-known that the $n$-th row of Pascal's triangle contains exactly $2^s$ odd binomial coefficients. We will show that the sum of these odd binomial coefficients is divisible by $2^s$. This seems to be a bit trickier to prove than the simplicity of its statement might suggest. We establish it, together with some generalisations, in Section \ref{app}.

It is not unusual for a constant-term sequence or a sequence defined by means of hypergeometric sums that the $p$-adic valuation of the $n$-th term should exhibit a relationship with the expansion of $n$ in base $p$. One quite general result to this effect is due to Delaygue \cite{D}. It applies to the counting sequences $W^\star_d(n)$ of the Abelian squares, and in terms of the counting sequences $W_d(n)$ of the high-dimensional lattice walks it tells us that $w_d(n) \ge 2s_2(n)$ if $d$ is even. Our Theorem shows that the latter inequality becomes an identity when $\nu_2(d) = 1$; while $w_d(n)$ is in fact much bigger when $\nu_2(d) \ge 2$.

The rest of the paper is structured as follows: We resolve parts (a) and (b) of the Theorem in Section \ref{init}. Sections \ref{abac-i}--\ref{abac-ii} introduce the aforementioned general technique, with Section \ref{app} collecting some applications. We give our first proof of parts (c) and (d) of the Theorem in Section \ref{proof-i}. The second proof comes next, in Section \ref{proof-ii}. Finally, in Section \ref{further}, we circle back to Bhattacharya's conjecture which caused us to embark on this entire train of thought and we discuss the question of whether our methods could perhaps shed some light on it as well.

\section{Initial Observations} \label{init}

We will need Kummer's theorem. Given a prime number $p$ and nonnegative integers $x_1$, $x_2$, $\ldots$, $x_k$ with sum $n$, it states that the exponent of $p$ in $\binom{n}{x_1, x_2, \ldots, x_k}$ equals $\frac{1}{p - 1}(s_p(x_1) + s_p(x_2) + \cdots + s_p(x_k) - s_p(n))$. Here, $s_p(x)$ is the sum of the digits of $x$ in base $p$.

The quantity $\frac{1}{p - 1}(s_p(x_1) + s_p(x_2) + \cdots + s_p(x_k) - s_p(n))$ equals the total number of carries which occur when $x_1$, $x_2$, $\ldots$, $x_k$ are added together in base $p$. The theorem is often stated in terms of this total number of carries instead of digit sums. One useful corollary is that $p$ does not divide $\binom{n}{x_1, x_2, \ldots, x_k}$ if and only if $x_1$, $x_2$, $\ldots$, $x_k$ form a carry-free partitioning of $n$.

We now specialise the theorem to $p = 2$. In this setting, the carry-free partitionings of $n$ are straightforward to describe explicitly. Indeed, split $n$ into $s_2(n)$ pairwise distinct powers of~$2$. Each one of these powers of $2$ must go to one piece in the partitioning. So the number of (ordered) carry-free partitionings of $n$ with $k$ pieces equals $k^{s_2(n)}$.

We continue with one simple but important observation. Recall that $W_d(n) = \binom{2n}{n}W^\star_d(n)$. By Kummer's theorem, the exponent of $2$ in $\binom{2n}{n}$ equals $s_2(n)$. Hence, writing $w^\star_d(n)$ for the exponent of $2$ in $W^\star_d(n)$, we arrive at \[w_d(n) = s_2(n) + w^\star_d(n).\]

This allows us to focus on $w^\star_d(n)$ in place of $w_d(n)$. Part (a) of the Theorem does not pose much difficulty anymore:

\begin{myproof}[Proof of part (a) of the Theorem] Suppose that $d$ is odd. We must show that $w^\star_d(n) = 0$. We use the formula ($\mathbf{W}^\star$). By the preceding discussion, the $(x_1, x_2, \ldots, x_d)$-summand in it is odd if and only if $x_1$, $x_2$, $\ldots$, $x_d$ form a carry-free partitioning of $n$; and there exist a total of $d^{s_2(n)}$ such partitionings. Since the latter count is odd, so must be $W^\star_d(n)$ as well. \end{myproof}

Suppose, from now on, that $d$ is even and $d = 2g$. For part (b), first of all we must work out a more helpful formula for $W^\star_d(n)$. Let $y_i = x_{2i - 1} + x_{2i}$ for all $i = 1$, $2$, $\ldots$, $g$. Then \[\binom{n}{x_1, x_2, \ldots, x_d} = \binom{n}{y_1, y_2, \ldots, y_g} \binom{y_1}{x_1, x_2} \binom{y_2}{x_3, x_4} \cdots \binom{y_g}{x_{d - 1}, x_d}.\]

Furthermore, by Vandermonde's identity, $\binom{y}{0}^2 + \binom{y}{1}^2 + \cdots + \binom{y}{y}^2 = \binom{2y}{y}$ for all nonnegative integers $y$; and, in particular, for all $y = y_1$, $y_2$, $\ldots$, $y_g$. We apply both of these observations to the formula ($\mathbf{W}^\star$), in order, and we arrive at \[W^\star_d(n) = \sum_{y_1 + y_2 + \cdots + y_g = n} \binom{n}{y_1, y_2, \ldots, y_g}^2 \binom{2y_1}{y_1} \binom{2y_2}{y_2} \cdots \binom{2y_g}{y_g}. \tag{$\mathbf{W}^\mathsection$}\]

Once the ``right'' formula has been found, part (b) yields to more or less the same method as part (a), with minimal adaptations.

\begin{myproof}[Proof of part (b) of the Theorem] Suppose that $\nu_2(d) = 1$ and $g$ is odd. We must show that $w^\star_d(n) = s_2(n)$. We use the formula ($\mathbf{W}^\mathsection$). By Kummer's theorem, the exponent of $2$ in its $(y_1, y_2, \ldots, y_g)$-summand equals $3(s_2(y_1) + s_2(y_2) + \cdots + s_2(y_g)) - 2s_2(n)$. Since $s_2(y_1) + s_2(y_2) + \cdots + s_2(y_g) \ge s_2(n)$, the latter quantity is always at least $s_2(n)$; and equality is attained if and only if $y_1$, $y_2$, $\ldots$, $y_g$ form a carry-free partitioning of $n$. To show that the smallest exponent of $2$ in an individual summand equals the exponent of $2$ in the whole sum, we must verify that the minimum is attained an odd number of times. This is indeed true, as the number of carry-free partitionings of $n$ with $g$ pieces equals~$g^{s_2(n)}$. \end{myproof}

The same argument demonstrates also that $w^\star_d(n) \ge s_2(n)$ and $w_d(n) \ge 2s_2(n)$ whenever $d$ is even.

Parts (c) and (d) of the Theorem, where $\nu_2(d) \ge 2$, are substantially more difficult. We begin our study of them in the next section. The discussion will initially focus on the special case when $d = 4$. However, soon we are going to see that the tools we have developed for this narrow purpose do in fact generalise to a much wider range of problems.

\section{Sums over Abaci I} \label{abac-i}

Say we are given a sum and we wish to show that it is divisible by a certain prime power. Our central example, motivating most of the discussion, will be the sum \[W^\star_4(n) = \sum_{x + y = n} \binom{n}{x, y}^2 \binom{2x}{x} \binom{2y}{y} \tag{$\mathbf{W}^\mathsection_4$}\] obtained when ($\mathbf{W}^\mathsection$) is specialised with $d = 4$. For convenience, let $s = s_2(n)$. We wish to show that this sum (which equals the $n$-th Domb number) is divisible by $2^{2s}$.

Observe that each summand of ($\mathbf{W}^\mathsection_4$) is already divisible by a pretty high power of $2$. Specifically, the $(x, y)$-summand $E(x, y)$ is divisible by $2$ with exponent $\omega(x, y) = 3(s_2(x) + s_2(y)) - 2s_2(n) \ge s$ as in the proof of part (b) of the Theorem.

The summands with a big $\omega$ make our task easier, while the summands with a small $\omega$ make it harder. It is natural now to imagine all summands of ($\mathbf{W}^\mathsection_4$) sorted into classes depending on how much difficulty they contribute to the problem. We ought to expect the summands where $\omega$ is minimal to make the most trouble for us, and so it is with them that we begin our analysis.

Consider the sub-sum $S$ of ($\mathbf{W}^\mathsection_4$) obtained when all such summands are split off from the rest. The necessary and sufficient condition for $\omega$ to attain its minimum is that, when we add $x$ and $y$ together, the computation is carry-free in binary. Hence, $S$ is the sum of $E(x, y)$ over all carry-free partitionings of $n$ into two nonnegative integers $x$ and $y$.

Our intuition is that the exponent of $2$ in the prime factorisation of $S$ ought to give us a measure of how much trouble the toughest summands of ($\mathbf{W}^\mathsection_4$) are going to make for us. This exponent is straightforward to determine experimentally. We examine its values over the first few positive integers $n$, and we find that it seems to never fall below $2s$. This is quite the lucky turn of events for us, and we are motivated by it to search for a generalisation of the same effect to the remaining ``difficulty classes''. Eventually, we arrive at the following notion:

We define an \emph{abacus} to be a two-row $\{0, 1\}$-matrix. Given an ordered pair of nonnegative integers $(x, y)$, we associate with it the abacus obtained by putting the binary expansions of $x$ and $y$ in the first and second row, respectively. When these binary expansions are not of the same length, we pad the shorter one with meaningless leading zeroes as needed. Note that the abacus of $x$ and $y$ is precisely how we would arrange their binary expansions so as compute $x + y$ by hand, with pen and paper.

We label each column in an abacus with the letter $\sfZ$ when it contains two $0$s; $\sfO$ when it contains two $1$s; and $\sfM$ when it is mixed, i.e., when the two binary digits in it are distinct. The \emph{type} of an abacus is the word over the alphabet $\{\sfZ, \sfO, \sfM\}$ obtained when we read off these labels, in order from left to right.

The value of $x + y$ is constant over all abaci of the same type $T$; we denote it by $\sigma(T)$. Let also $\alpha_I(T)$ be the number of occurrences of the letter $I$ in $T$. Then each type $T$ is associated with a total of $2^{\alpha_\sfM(T)}$ abaci. We write $\varepsilon$ for the empty type; the unique abacus of it is the one with zero columns, when $x = y = 0$. We also write $\cf(n)$ for the type obtained from the binary expansion of $n$ by replacing each $0$ with a $\sfZ$ and each $1$ with an $\sfM$. A partitioning of $n$ into two nonnegative integer pieces $x$ and $y$ is carry-free if and only if the abacus associated with it is of type $\cf(n)$.

Given a two-variable function $F$ and an abacus type $T$, we write $\sum_T F(x, y)$ for the sum of $F$ over all ordered pairs of nonnegative integers $x$ and $y$ such that the corresponding abacus is of type $T$. Every sum of the form $\sum_{x + y = n} F(x, y)$ can be split into such sub-sums over abacus types.

We return now to our motivational example. Using our new notation, we can rephrase our experimental observations as follows: ``$2^{2s}$ seems to divide $\sum_{\cf(n)} E(x, y)$''. Our task will be accomplished if we manage to show that, in fact, $2^{2s}$ divides $\sum_T E(x, y)$ for all abacus types $T$ with $n = \sigma(T)$.

We begin by taking a closer look at the structure of $E(x, y)$. For our purposes, the factorials will turn out to be a lot more convenient as building blocks than the binomial coefficients, and so we rewrite $E(x, y)$ as \[n!^2 \cdot (x!y!)^{-4} \cdot (2x)!(2y)!.\]

Since $s_2(x) + s_2(y)$ remains constant over each abacus type $T$, so does $\omega(x, y)$ as well. Hence, it makes sense to factor the corresponding power of $2$ out of all summands in $\sum_T E(x, y)$. Given a positive integer $z$, we denote its odd part by $\Theta_2(z)$. Then $\Theta_2(n!)^2$ can be factored out as well, and we are left with \[\big( \Theta_2(x!)\Theta_2(y!) \big)^{-4} \cdot \Theta_2((2x)!)\Theta_2((2y)!).\]

We continue by embedding this ``residue'' into one wider family of similarly-structured functions. It is this more general family that we are going to draw on for applications.

Fix a symmetric integer-coefficient two-variable polynomial $P$. Fix also a vector $\bfe = (e_0,\allowbreak e_1, \ldots, e_k)$ with nonnegative integer components. Suppose, from now on, that \[F(x, y) = P(x, y) \cdot \prod_{i = 0}^k \big( \Theta_2((2^ix)!) \Theta_2((2^iy)!) \big)^{e_i}.\]

(Notice that the aforementioned residue is not, strictly speaking, a member of this family. Formally, in order to obtain it we must set $P = 1$ and $\bfe = (-4, 1)$; however, one of the components of $\bfe$ will then be negative. This should not bother us too much. Since each odd number is invertible modulo all powers of $2$, in some sense we can replace all negative components of $\bfe$ with positive ones. We sort out the details near the end of this section.)

We introduce also the notation \[\mathfrak{S}^\bfe_T(P) = \sum_T F(x, y).\]

We proceed now to derive one system of reduction formulas for sums of this form. Roughly speaking, our goal will be to make $T$ smaller at the cost of making $P$ more complicated. By means of a series of such reductions, eventually we are going to reach a sum where $T$ is the empty type whereas $P$ has ballooned in complexity. Of course, in the latter setting $\mathfrak{S}^\bfe_\varepsilon(P) = P(0, 0)$.

First we must take care of some preliminaries. Let $\psi_i(x) = \Theta_2((2^ix)!)$.

\begin{lemma} \label{fact} For each nonnegative integer $i$, there exists an integer-coefficient polynomial $\Psi_i$ such that $\psi_i(2x + 1) = \psi_{i + 1}(x)\Psi_i(x)$. \end{lemma} 

\begin{myproof} By \[\frac{\psi_i(2x + 1)}{\psi_{i + 1}(x)} = \prod_{j = 1}^{2^i} \Theta_2(2^{i + 1}x + j) = \prod_{j = 1}^{2^i} \big( 2^{i + 1 - \nu_2(j)}x + \Theta_2(j) \big);\] each factor on the right-hand side is a linear function of $x$ with integer coefficients. \end{myproof}

For convenience, set \[\Psi^\bfe(x) = \prod_{i = 0}^k \Psi_i(x)^{e_i}.\]

We define \begin{gather*} Q_\sfZ(x, y) = P(2x, 2y)\\ Q_\sfO(x, y) = \Psi^\bfe(x)\Psi^\bfe(y)P(2x + 1, 2y + 1)\\ Q_\sfM(x, y) = \frac{1}{2}\big( \Psi^\bfe(x)P(2x + 1, 2y) + \Psi^\bfe(y)P(2x, 2y + 1) \big). \end{gather*}

We are also going to need the notation $\shift \bfe = (0, e_0, e_1, \ldots, e_k)$.

\begin{lemma} \label{red} For all $I$, it holds that $Q_I$ is a symmetric integer-coefficient polynomial. Furthermore, if $I \in \{\sfZ, \sfO\}$, then \[\mathfrak{S}^\bfe_{TI}(P) = \mathfrak{S}^{\shift \bfe}_T(Q_I);\] and, if $I = \sfM$, then \[\mathfrak{S}^\bfe_{T \sfM}(P) = 2 \cdot \mathfrak{S}^{\shift \bfe}_T(Q_\sfM).\] \end{lemma} 

\begin{myproof} The identities are immediate from the definitions. The fact that all coefficients of $Q_I$ are integers is immediate, too, if $I \in \{\sfZ, \sfO\}$. Otherwise, if $I = \sfM$, it suffices to show that all coefficients of $2Q_\sfM$ are even integers. The proof of Lemma \ref{fact} makes it obvious that $\Psi^\bfe(x) \equiv \Psi^\bfe(y) \equiv 1 \pmod 2$. Hence, $2Q_\sfM = \Psi^\bfe(x)P(2x + 1, 2y) + \Psi^\bfe(y)P(2x, 2y + 1) \equiv P(1, 0) + P(0, 1) \pmod 2$. The right-hand side works out to an even integer because $P$ is symmetric. \end{myproof}

Right away, these reductions give us the following:

\begin{lemma} \label{div} The exponent of $2$ in $\mathfrak{S}^\bfe_T(P)$ is at least $\alpha_\sfM(T)$. \end{lemma} 

\begin{myproof} By iterating the reductions of Lemma \ref{red}, we can find an integer vector $\bff$ and a symmetric integer-coefficient two-variable polynomial $R$ such that \[\mathfrak{S}^\bfe_T(P) = 2^{\alpha_\sfM(T)} \cdot \mathfrak{S}^{\bff}_\varepsilon(R).\]

Since the factor $\mathfrak{S}^{\bff}_\varepsilon(R)$ on the right-hand side is an integer, the result follows. \end{myproof}

Lemma \ref{div} is already strong enough to settle our motivational example -- i.e., to show that $2^{2s}$ divides the right-hand side of ($\mathbf{W}^\mathsection_4$). We spell out the details in the next section. However, for the proof of our main Theorem, we require also some way to pin down the cases of equality. Or, in other words, the cases when $\alpha_\sfM(T) = \nu_2(\mathfrak{S}^\bfe_T(P))$ exactly.

Below, we will see that the equality cases of Lemma \ref{div} are described by a certain finite automaton. The gist of the proof is that it suffices to keep track of just a tiny bit of information about $P$, over the course of our reductions, so as to be able to tell at the end whether $\mathfrak{S}^{\bff}_\varepsilon(R)$ is even or not.

We go on to introduce special notation for the tiny pieces of information we will be keeping track of. Let $\widehat{\bfe} = (e_0, e_1, e_2 + e_3 + \cdots + e_k)$. We call this the \emph{fold} of $\bfe$. Let also $\widehat{P} = P \bmod 2$ and $\widehat{Q_I} = Q_I \bmod 2$ for all $I$. (Given a positive integer $t$ and an integer-coefficient polynomial $\Phi$, we write $\Phi \bmod t$ for the remainder of $\Phi$ modulo $t$. It is defined as the unique polynomial $\widehat{\Phi}$ with coefficients in $\{0, 1, \ldots, t - 1\}$ such that $\Phi \equiv \widehat{\Phi} \pmod t$.)

\begin{lemma} \label{step} If we know $\widehat{\bfe}$ and $\widehat{P}$, we can determine $\widehat{Q_I}$ for all $I$ without any additional information about $\bfe$ and $P$. Furthermore, for all $I$, it holds that $\deg \widehat{Q_I} \le \max\{1, \deg \widehat{P}\}$. \end{lemma}

\begin{myproof} The claim is obvious if $I \in \{\sfZ, \sfO\}$. We assume from now on that $I = \sfM$. Then $2Q_\sfM$ is an integer-coefficient polynomial. We must show that $2Q_\sfM \bmod 4$ is uniquely determined by $\widehat{\bfe}$ and $\widehat{P}$, and we must also estimate its degree.

Let $P = 2P_\divideontimes + \widehat{P}$. Then $\Psi^\bfe(x)P_\divideontimes(2x + 1, 2y) + \Psi^\bfe(y)P_\divideontimes(2x, 2y + 1) \equiv 0 \pmod 2$ by the same argument as in the proof of Lemma \ref{red}. So $2Q_\sfM = \Psi^\bfe(x)P(2x + 1, 2y) + \Psi^\bfe(y)P(2x, 2y + 1) \equiv \Psi^\bfe(x)\widehat{P}(2x + 1, 2y) + \Psi^\bfe(y)\widehat{P}(2x, 2y + 1) \pmod 4$.

Using the explicit form of $\Psi_i(x)$ in the proof of Lemma \ref{fact}, it is straightforward to check that $\Psi_0(x) \equiv \Psi_1(x) \equiv 2x + 1 \pmod 4$ and $\Psi_i(x) \equiv 2x + 3 \pmod 4$ when $i \ge 2$. So $\Psi^\bfe(x) \equiv \Psi^{\widehat{\bfe}}(x) \pmod 4$ and $\Psi^\bfe(y) \equiv \Psi^{\widehat{\bfe}}(y) \pmod 4$. Of course, given $\widehat{P}$, we can also find $\widehat{P}(2x + 1, 2y) \bmod 4$ and $\widehat{P}(2x, 2y + 1) \bmod 4$. Hence, it is indeed possible to determine $\widehat{Q_\sfM}$ based solely on $\widehat{\bfe}$ and $\widehat{P}$.

Observe next that both of $\Psi^\bfe(x)$ and $\widehat{P}(2x + 1, 2y)$ are constant modulo $2$; explicitly, $\Psi^\bfe(x) \equiv 1 \pmod 2$ and $\widehat{P}(2x + 1, 2y) \equiv P(1, 0) \pmod 2$. So we can find two integer-coefficient polynomials $A$ and $B$ as well as two integer constants $a$ and $b$ such that $\Psi^\bfe(x) = 2A(x) + a$ and $\widehat{P}(2x + 1, 2y) = 2B(x, y) + b$. Then $\Psi^\bfe(x)\widehat{P}(2x + 1, 2y) \equiv 2bA(x) + 2aB(x, y) + ab \pmod 4$; and analogous reasoning applies also to $\Psi^\bfe(y)\widehat{P}(2x, 2y + 1)$. We conclude that $\deg \widehat{Q_\sfM} \le \max\{\deg {(\Psi^\bfe(x) \bmod 4)},\allowbreak \deg {(\Psi^\bfe(y) \bmod 4)}, \deg \widehat{P}\} = \max \{1, \deg \widehat{P}\}$, as desired. \end{myproof}

We assume, for convenience, that all of our automata go through their input from right to left. This is the same order as the one in which our reductions disassemble the abacus types.

\begin{lemma} \label{auto} Given $\bfe$ and $P$, there exists a finite automaton $\mathcal{F}$ over the alphabet $\{\sfZ, \sfO, \sfM\}$ such that the exponent of $2$ in $\mathfrak{S}^\bfe_T(P)$ equals $\alpha_\sfM(T)$ if and only if $\mathcal{F}$ accepts $T$. \end{lemma} 

\begin{myproof} In the setting of the proof of Lemma \ref{div}, we must work out whether $\mathfrak{S}^{\bff}_\varepsilon(R) = R(0, 0)$ is even or not. For this purpose, it suffices to determine $\widehat{R} = R \bmod 2$.

Over the course of our reductions, the original sum $\mathfrak{S}^\bfe_T(P)$ is transformed into a succession of intermediate sums of the form $\mathfrak{S}^{\bfe_i}_{T_i}(P_i)$ with $\bfe_i = \shift^i \bfe$ and $T_i$ being the prefix of $T$ obtained when its rightmost $i$ letters are deleted. Suppose we manage to show that, for each one of $\widehat{\bfe_i}$ and $\widehat{P_i}$, there are only finitely many values which it could potentially take on during the reduction process. Then Lemma \ref{step} will guarantee the existence of a finite automaton $\mathcal{F}$ which, based on $T$, calculates $\widehat{R}$. (Formally, each state of $\mathcal{F}$ will be labelled with an ordered pair consisting of a fold and a remainder.)

This is clear in the case of $\widehat{\bfe_i}$. For $\widehat{P_i}$, it follows because all of its coefficients are in $\{0, 1\}$ and its degree is bounded from above as per Lemma \ref{step}. \end{myproof}

What remains is to deal away with the non-negativity restriction we imposed on $\bfe$. We define the exponent of $2$ in the rational number $u/v$ by $\nu_2(u/v) = \nu_2(u) - \nu_2(v)$.

\begin{lemma} \label{neg} Both Lemmas \ref{div} and \ref{auto} remain true when the components of $\bfe$ are allowed to be arbitrary (possibly negative) integers. \end{lemma} 

\begin{myproof} By induction on $i \ge 1$, we get that $z^{2^i} \equiv 1 \pmod{2^{i + 2}}$ for all odd $z$. Let $\bff_i = (f_{i, 0}, f_{i, 1}, \ldots, f_{i, k})$ with $f_{i, j} = e_j \bmod{2^i}$ for all $j$. Then $\mathfrak{S}^\bfe_T(P) \equiv \mathfrak{S}^{\bff_i}_T(P) \pmod{2^{i + 2}}$. Since the original Lemma \ref{div} applies to the right-hand side, its generalisation with $\bfe$ unrestricted follows when $i \ge \alpha_\sfM(T)$.

Consider next the finite automaton $\mathcal{F}_i$ obtained when we apply the construction of the original Lemma~\ref{auto} to $\bff_i$ and $P$. Observe that, in the proof of Lemma \ref{step}, the only way in which $\bfe$ enters into the calculation of $\widehat{Q_I}$ is via the remainders $\Psi^\bfe(x) \bmod 4$ and $\Psi^\bfe(y) \bmod 4$. However, $\Psi_j(x)^2 \equiv \Psi_j(y)^2 \equiv 1 \pmod 4$ for all $j$. Thus there exists a finite automaton $\mathcal{F}$ such that $\mathcal{F}$ and $\mathcal{F}_i$ are isomorphic for all sufficiently large $i$. This finite automaton will fit the bill, in the setting of the generalised Lemma \ref{auto} with an unrestricted $\bfe$, regardless of how big we make $\alpha_\sfM(T)$. \end{myproof}

\section{Applications} \label{app}

We are all set now to browse through some applications of the apparatus developed in the previous section. We begin with the one which was advertised in the introduction:

\begin{proposition} \label{odd} Let $n$ be a positive integer. Consider the sum of all odd binomial coefficients in the $n$-th row of Pascal's triangle. The exponent of $2$ in this sum is at least $s_2(n)$. Equality is attained if and only if the binary expansion of $n$ does not contain two adjacent $1$s. \end{proposition}

\begin{table}[ht] \centering \footnotesize \begin{tabular}{|c|c|c|}
from & with & to\\
\hline
\multirow{2}{*}{$1$} & $0$ & $1$\\
\cline{2-3}
& $1$ & $x + y + 1$\\
\hline
\multirow{2}{*}{$x + y + 1$} & $0$ & $1$\\
\cline{2-3}
& $1$ & $0$\\
\hline
$0$ & $0$, $1$ & $0$\\
\hline
\end{tabular} \caption{} \label{sub-11} \end{table}

\begin{myproof} By Lemma \ref{neg}, as our sum equals $\Theta_2(n!) \cdot \mathfrak{S}^{(-1)}_{\cf(n)}(1)$. The associated finite automaton is shown in Table \ref{sub-11}. We use the alphabet $\{0, 1\}$ instead of $\{\sfZ, \sfO, \sfM\}$ because it is only the abacus type $\cf(n)$ which matters for our purposes, and the binary expansion of $n$ maps onto it by means of the substitutions $0 \to \sfZ$ and $1 \to \sfM$. The general construction prescribes that each state must be labelled with a fold and a remainder; but, in this instance, states with the same remainder behave identically. \end{myproof}

Our considerations in the previous section were motivated throughout by the problem of showing that the $n$-th Domb number is divisible by $2^{2s_2(n)}$. We resolve this problem as follows:

\begin{proposition} \label{dn} The exponent of $2$ in the $n$-th Domb number is at least $2s_2(n)$. Equality is attained if and only if the binary expansion of $n$ does not contain two adjacent $1$s. \end{proposition}

\begin{myproof} Consider any abacus type $T$ with $n = \sigma(T)$. For convenience, set $\alpha_1 = \alpha_\sfM(T)$ and $\alpha_2 = \alpha_\sfO(T)$. So $s_2(x) + s_2(y) = \alpha_1 + 2\alpha_2$ over $T$; and $s_2(n) \le \alpha_1 + \alpha_2$. Define $\omega$ and $E$ as in Section \ref{abac-i}. We get that $\omega$ is constant over $T$. Denote its value by $\omega(T)$. Then $\omega(T) = 3(s_2(x) + s_2(y)) - 2s_2(n) \ge \alpha_1 + 4\alpha_2$.

Partition ($\mathbf{W}^\mathsection_4$) into sub-sums of the form $\sum_T E(x, y) = 2^{\omega(T)} \cdot \Theta_2(n!)^2 \cdot \mathfrak{S}^{(-4, 1)}_T(1)$. By Lemma~\ref{neg}, the last factor on the right-hand side is divisible by~$2^{\alpha_1}$. Hence, the exponent of $2$ in $\sum_T E(x, y)$ is always at least $\omega(T) + \alpha_1 \ge 2\alpha_1 + 4\alpha_2 \ge 2s_2(n) + 2\alpha_2$. The latter quantity will be strictly greater than $2s_2(n)$ unless $\alpha_2 = 0$; i.e., unless $T = \cf(n)$. The associated finite automaton is once again as shown in Table~\ref{sub-11}. \end{myproof}

In response to the author's queries as to the originality of Proposition \ref{odd}, it was noted \cite{P} that the divisibility can be obtained as a corollary of one earlier result -- specifically, the first part of Lemma 12 in Calkin's \cite{C}. It states that, for each positive integer $e$, the exponent of $2$ in $\binom{n}{0}^e + \binom{n}{1}^e + \cdots + \binom{n}{n}^e$ is at least $s_2(n)$. (Though it does not address the necessary and sufficient conditions for equality.) The derivation of the divisibility part of Proposition \ref{odd} therefrom is by setting $e = 2^k + 1$ with $k$ sufficiently large. We show now that our methods can be used to tackle the more general result as well:

\begin{proposition} \label{pow} Let $e$ and $n$ be positive integers with $e \ge 2$. Consider the sum of the $e$-th powers of all binomial coefficients in the $n$-th row of Pascal's triangle. When $e$ is even, the exponent of $2$ in this sum equals $s_2(n)$. Otherwise, when $e$ is odd, the same exponent is at least $s_2(n)$; and equality is attained if and only if the binary expansion of $n$ does not contain two adjacent~$1$s. \end{proposition} 

\begin{myproof} Consider any abacus type $T$ with $n = \sigma(T)$. Define $\alpha_1$ and $\alpha_2$ as in the proof of Proposition~\ref{dn}. By Kummer's theorem, the exponent of $2$ in $\binom{n}{x, y}^e$ equals $e(s_2(x) + s_2(y) - s_2(n))$, and so it is constant over $T$. Denote its value by $\xi(T)$. Then $\xi(T) \ge e\alpha_2$.

Partition the full sum into sub-sums of the form $\sum_T \binom{n}{x, y}^e = 2^{\xi(T)} \cdot \Theta_2(n!)^e \cdot \mathfrak{S}^{(-e)}_T(1)$. By Lemma \ref{neg}, the last factor on the right-hand side is divisible by~$2^{\alpha_1}$. Hence, the exponent of $2$ in $\sum_T \binom{n}{x, y}^e$ is always at least $\xi(T) + \alpha_1 \ge \alpha_1 + e\alpha_2 \ge s_2(n) + (e - 1)\alpha_2$. The rest of the argument proceeds along the same lines as in the proof of Proposition \ref{dn}. When $e$ is odd, we obtain the same finite automaton as before. Otherwise, when $e$ is even, after minimisation we obtain a finite automaton with a single state which all transition arrows loop back to. \end{myproof}

The work \cite{OS} is about a generalisation of the Domb numbers given by \[D_{a, b, c}(n) = \sum_{x + y = n} \binom{n}{x, y}^a \binom{2x}{x}^b \binom{2y}{y}^c.\]

Our methods are applicable when $b = c$, and we arrive at the following result which generalises both Propositions \ref{dn} and \ref{pow}:

\begin{proposition} \label{gdn} Let $a$ and $b$ be nonnegative integers with $a + b \ge 2$. Let also $n$ be a positive integer, and consider the generalised Domb number $D_{a, b, b}(n)$. When $a + b$ is even, the exponent of $2$ in it equals $(b + 1)s_2(n)$. Otherwise, when $a + b$ is odd, the same exponent is at least $(b + 1)s_2(n)$; and equality is attained if and only if the binary expansion of $n$ does not contain two adjacent~$1$s. \end{proposition}

We omit the proof as it proceeds along the same lines as those of Propositions \ref{dn} and \ref{pow}. Observe that, if $a + b \le 1$, then the generalised Domb number $D_{a, b, b}(n)$ can be calculated explicitly. The identities $D_{0, 0, 0}(n) = n + 1$ and $D_{1, 0, 0}(n) = 2^n$ are clear; whereas it is well-known that $D_{0, 1, 1}(n) = 2^{2n}$.

We now skip ahead a bit. Propositions \ref{mult}--\ref{prime} below require the fully general form of our apparatus developed in the next section. For this reason, we limit ourselves to quick sketches in place of complete proofs.

\begin{proposition} \label{mult} Let $k$ and $n$ be positive integers. Consider the sum of all odd $2k$-nomial coefficients with $n$ at the top. The exponent of $2$ in this sum is at least $s_2(n) + \nu_2(k)$. When $k$ is odd, equality is attained if and only if the binary expansion of $n$ does not contain two adjacent~$1$s. Otherwise, when $k$ is even, equality is attained if and only if, in the binary expansion of $n$, all of the $1$s come before all of the $0$s. \end{proposition} 

This is a generalisation of Proposition \ref{odd} to multinomial coefficients which bears a strong resemblance to parts (c) and (d) of our main Theorem. The same method works as in Section~\ref{proof-i} -- specifically, that of halting the reductions not at $\varepsilon$ but at a single-letter type instead.

We can also generalise Proposition \ref{odd} in other directions, by considering different primes:

\begin{proposition} \label{tri} Let $n$ be a positive integer. Consider the sum of all trinomial coefficients with $n$ at the top which are not divisible by $3$. The exponent of $3$ in this sum is always at least the number of nonzero digits in the ternary expansion of $n$. Equality is attained if and only if this ternary expansion begins with $1$ but does not contain the subwords $11$ and $02$. \end{proposition} 

Here, matters are complicated by the fact that the relevant trinomial coefficients do not all belong to the same type. One way to get around this difficulty is to modify the notion of an abacus type so that each letter encodes the sum of the digits in the corresponding abacus column rather than the multiset formed by these digits. The discussion in the next section continues to apply when this modification is implemented with $p = 3$ and $r = 3$. The trinomial coefficients which make up our sum will then indeed be captured by a single one of these new abacus types.

(It is worth noting also that the associated finite automaton supplied by the ``standard'' construction can be minimised by disregarding all folds and merging together the states whose remainders differ solely by a nonzero constant factor.)

The directions of generalisation suggested by Propositions \ref{mult} and \ref{tri} can be combined, too:

\begin{proposition} \label{prime} Let $p$ be a prime number. Let also $k$ and $n$ be positive integers. Consider the sum of all $kp$-nomial coefficients with $n$ at the top which are not divisible by $p$. The exponent of $p$ in this sum is always at least $s^\star_p(n) + \nu_p(k)$, where $s^\star_p(n)$ is the number of nonzero digits in the expansion of $n$ to base $p$. Furthermore, with $k$ and $p$ fixed, there exists a finite automaton $\mathcal{F}$ such that equality is attained if and only if $\mathcal{F}$ accepts the expansion of $n$ to base $p$. \end{proposition}

The same difficulties arise here as with Proposition \ref{tri}, and we modify the basic notion of an abacus type in the same manner as before in order to resolve them. Following this, once again we apply the ``early halting'' method of Section \ref{proof-i}.

Observe that, for all $k$ and $p$, equality is attained whenever $n$ is a power of $p$. Thus the language of $\mathcal{F}$ is never empty. We do not, however, attempt an explicit description of the equality cases in this setting. It would be interesting to know whether such a description is in fact possible. To arrive at it, one would need to analyse the infinite family of finite automata $\mathcal{F}$ obtained when $k$ and $p$ vary.

\section{Sums over Abaci II} \label{abac-ii}

Here, we generalise the material of Section \ref{abac-i}. Many of the terms and notations will be revised so that the originals become special cases.

Fix a prime number $p$ and a positive integer $r$. We define an \emph{abacus} to be a matrix with $r$ rows where each entry is a digit in base $p$. (We had $p = 2$ and $r = 2$ in Section \ref{abac-i}.) Given a vector $\bfx = (x_1, x_2, \ldots, x_r)$ with nonnegative integer components, the abacus of $\bfx$ is obtained by putting the expansion of $x_i$ to base $p$ into row $i$. As before, we pad these expansions with meaningless leading zeroes as needed so as to ensure that all of them are of the same length.

Let $\mathcal{A}$ be a finite alphabet whose letters are in a one-to-one correspondence with the multisets of size $r$ consisting of digits in base $p$. We label each abacus column with the letter of $\mathcal{A}$ which describes its contents. The \emph{type} of an abacus is the word over $\mathcal{A}$ obtained by reading off these letters, in order from left to right.

Let $F$ be a function of the variables $x_1$, $x_2$, $\ldots$, $x_r$; we use $F(\bfx)$ as shorthand for $F(x_1, x_2,\allowbreak \ldots, x_r)$. Given a type $T$, we write $\sum_T F(\bfx)$ for the sum of $F$ over all $\bfx$ such that the abacus of $\bfx$ is of type $T$. As in Section \ref{abac-i}, we are going to focus on sums of this kind with $F$ of one particular form, which we proceed now to specify.

Given a positive integer $x$, we define $\Theta_p(x)$ by $x = p^{\nu_p(x)}\Theta_p(x)$. Let also $\psi_i(x) = \Theta_p((p^ix)!)$.

Fix a symmetric integer-coefficient $r$-variable polynomial $P$ as well as a vector $\bfe = (e_0, e_1,\allowbreak \ldots, e_k)$ with nonnegative integer components. (As in Section \ref{abac-i}, the non-negativity restriction will eventually be lifted.) Then set \[F(\bfx) = P(\bfx) \cdot \prod_{i = 0}^k \big( \psi_i(x_1)\psi_i(x_2) \cdots \psi_i(x_r) \big)^{e_i}\] and \[\mathfrak{S}^\bfe_T(P) = \sum_T F(\bfx).\]

The next result generalises Lemma \ref{fact}.

\begin{lemma} \label{gfact} For each nonnegative integer $i$ and each digit $a$ in base $p$, there exists an integer-coefficient polynomial $\Psi_{i, a}$ such that $\psi_i(px + a) = \psi_{i + 1}(x)\Psi_{i, a}(x)$. \end{lemma} 

\begin{myproof} Explicitly, \[\Psi_{i, a}(x) = \prod_{j = 1}^{ap^i} \big( p^{i + 1 - \nu_p(j)}x + \Theta_p(j) \big),\] by the same calculations as in the proof of Lemma \ref{fact}. \end{myproof}

Given a vector $\bfa = (a_1, a_2, \ldots, a_r)$ all of whose components are digits in base $p$, we define \[\Psi^\bfe_\bfa(\bfx) = \prod_{i = 0}^k \big( \Psi_{i, a_1}(x_1)\Psi_{i, a_2}(x_2) \cdots \Psi_{i, a_r}(x_r) \big)^{e_i}.\]

We go on now to generalise the system of reduction formulas introduced in Section \ref{abac-i}. For each letter $I$ of $\mathcal{A}$, let $\varrho(I)$ be the set of all vectors whose components form a permutation of the elements of the multiset associated with $I$. We call $I$ \emph{special} when $p$ divides the size of $\varrho(I)$. Let $\lambda$ be the characteristic function of this property; i.e., $\lambda(I) = 1$ when $I$ is special and $\lambda(I) = 0$ otherwise. We also let $\lambda(T)$ be the sum of $\lambda(I)$ over all letters $I$ of $T$; or, equivalently, the sum of $\alpha_I(T)$ over all special letters $I$ of $\mathcal{A}$.

For each $I$, we set \[Q_I(\bfx) = p^{-\lambda(I)} \cdot \sum_{\bfa \in \varrho(I)} P(p\bfx + \bfa) \Psi^\bfe_\bfa(\bfx).\]

The analogues of Lemmas \ref{red} and \ref{div} are immediate:

\begin{lemma} \label{gred} For all $I$, it holds that $Q_I$ is a symmetric integer-coefficient polynomial with \[\mathfrak{S}^\bfe_{TI}(P) = p^{\lambda(I)} \cdot \mathfrak{S}^{\shift \bfe}_T(Q_I).\] \end{lemma} 

\begin{myproof} Just as in the proof of Lemma \ref{red}, the identity follows directly by the definitions. The fact that the coefficients of $Q_I$ are integers is straightforward, too, in the case when $I$ is non-special. Otherwise, if $I$ is special, we must verify that all coefficients of $pQ_I$ are integer multiples~of~$p$.

Since $P$ is symmetric, we get that $P(\bfa)$ remains constant when $\bfa$ varies over $\varrho(I)$. Denote its value by $P(I)$. In addition, the proof of Lemma \ref{gfact} shows that each $\Psi_{i, a}$ is constant modulo~$p$. Hence, there exists an integer constant $c$ such that $\Psi^\bfe_\bfa(\bfx) \equiv c \pmod p$ for all elements $\bfa$ of $\varrho(I)$. We conclude that $pQ_I = \sum_{\bfa \in \varrho(I)} P(p\bfx + \bfa) \Psi^\bfe_\bfa(\bfx) \equiv c \cdot |\varrho(I)| \cdot P(I) \pmod p$. However, $p$ divides the size of $\varrho(I)$ because $I$ is special, and the desired result follows. \end{myproof}

\begin{lemma} \label{gdiv} The exponent of $p$ in $\mathfrak{S}^\bfe_T(P)$ is at least $\lambda(T)$. \end{lemma} 

The proof is essentially identical to that of Lemma \ref{div}, and we omit it. What remains is to get a handle on the cases of equality in Lemma \ref{gdiv}. The gist of the argument will be exactly the same as before. However, the technical details will be slightly more complicated. For this reason, we sort out some of them in one preliminary lemma without a direct analogue in Section \ref{abac-i}.

\begin{lemma} \label{modp2} Suppose that $p$ is odd. Then the polynomials $\Psi_{i, a}$ satisfy the following congruences:

(a) $\Psi_{i, a}^{p^2 - p} \equiv 1 \pmod{p^2}$ for all $i \ge 0$ and $a$.

(b) $\Psi_{i, a} \equiv \Psi_{i + 2, a} \pmod{p^2}$ for all $i \ge 1$ and $a$. \end{lemma}

\begin{myproof} As in the proof of Lemma \ref{gfact}, $\Psi_{i, a}$ is the product of several linear functions with slopes divisible by $p$ but constant terms not divisible by $p$. Hence, for part (a) it suffices to show that $(u + vx)^{p^2 - p} \equiv 1 \pmod{p^2}$ whenever $u$ and $v$ are integers with $u \not \equiv v \equiv 0 \pmod p$. This is routine by expanding the left-hand side into $\sum_{j = 0}^{p^2 - p} \binom{p^2 - p}{j} u^{p^2 - p - j}v^jx^j$ and then, at $j = 0$, invoking Euler's theorem with $\varphi(p^2) = p^2 - p$.

We continue on to part (b). Given a positive integer $t$, let $\Gamma(t)$ be the product of all non-multiples of $p$ in the interval $[1; t]$. Clearly, $\Psi_{i + 1, a}(x)/\Psi_{i, a}(x)$ is the product of $p^{i + 2}x + j$ over all non-multiples $j$ of $p$ in the interval $[1; ap^{i + 1}]$. So $\Psi_{i + 1, a} \equiv \Gamma(ap^{i + 1}) \cdot \Psi_{i, a} \pmod{p^2}$. Similarly, $\Psi_{i + 2, a} \equiv \Gamma(ap^{i + 2}) \cdot \Psi_{i + 1, a} \pmod{p^2}$. Hence, for part (b) it suffices to demonstrate that $\Gamma(ap^{i + 1})\Gamma(ap^{i + 2}) \equiv 1 \pmod{p^2}$ when $i \ge 1$. For convenience, set $\gamma = \Gamma(p^2)$. Then $\Gamma(ap^{i + 1}) \equiv \gamma^{ap^{i - 1}} \pmod{p^2}$ and $\Gamma(ap^{i + 2}) \equiv \gamma^{ap^i} \pmod{p^2}$. However, $\gamma \equiv -1 \pmod{p^2}$ because all nonzero remainders modulo $p^2$ distinct from $\pm 1$ form several inverse pairs. \end{myproof}

We are ready to generalise Lemmas \ref{step} and \ref{auto}. The notion of a fold is left unchanged when $p = 2$. Otherwise, when $p$ is odd, we set $\widehat{\bfe} = (e_0, e_1 + e_3 + e_5 + \cdots, e_2 + e_4 + e_6 + \cdots)$. We define also $\widehat{P} = P \bmod p$ and $\widehat{Q_I} = Q_I \bmod p$ for all $I$.

\begin{lemma} \label{gstep} If we know $\widehat{\bfe}$ and $\widehat{P}$, we can determine $\widehat{Q_I}$ for all $I$ without any additional information about $\bfe$ and $P$. Furthermore, for all $I$, it holds that $\deg \widehat{Q_I} \le \max\{1, \deg \widehat{P}\}$. \end{lemma} 

\begin{myproof} The claim is straightforward when $I$ is non-special. Suppose, from now on, that it is in fact special. Then $pQ_I$ is an integer-coefficient polynomial. We must show that $pQ_I \bmod{p^2}$ is uniquely determined by $\widehat{\bfe}$ and $\widehat{P}$, and we must also estimate its degree.

By the same reasoning as in the proof of Lemma \ref{step}, we get that $pQ_I \equiv \sum_{\bfa \in \varrho(I)} \widehat{P}(p\bfx + \bfa) \Psi^\bfe_\bfa(\bfx) \pmod{p^2}$. Part (b) of Lemma \ref{modp2} tells us that $\Psi^\bfe_\bfa(x) \equiv \Psi^{\widehat{\bfe}}_\bfa(x) \pmod{p^2}$. Of course, given $\widehat{P}$, we can also find $\widehat{P}(p\bfx + \bfa) \bmod{p^2}$ for all elements $\bfa$ of $\varrho(I)$. Hence, it is indeed possible to determine $\widehat{Q_I}$ based solely on $\widehat{\bfe}$ and $\widehat{P}$.

Fix now an element $\bfa$ of $\varrho(I)$. Observe that both of $\Psi^\bfe_\bfa(\bfx)$ and $\widehat{P}(p\bfx + \bfa)$ are constant modulo~$p$. So we can find two integer-coefficient polynomials $U$ and $V$ as well as two integer constants $u$ and $v$ with $\Psi^\bfe_\bfa(\bfx) = pU(\bfx) + u$ and $\widehat{P}(p\bfx + \bfa) = pV(\bfx) + v$. Then $\widehat{P}(p\bfx + \bfa) \Psi^\bfe_\bfa(\bfx) \equiv pvU(\bfx) + puV(\bfx) + uv \pmod{p^2}$. Since the same reasoning applies to all $\bfa$, we conclude that $\deg \widehat{Q_I}$ cannot exceed the greatest one among $\deg {(\Psi^\bfe_\bfa(\bfx) \bmod{p^2})}$ and $\deg \widehat{P}(p\bfx + \bfa)$ over all elements $\bfa$ of $\varrho(I)$. The desired result follows as $\Psi^\bfe_\bfa(\bfx) \bmod{p^2}$ is always at most linear (by expanding the product in the proof of Lemma \ref{gfact}) and $\deg \widehat{P}(p\bfx + \bfa) = \deg \widehat{P}$. \end{myproof}

\begin{lemma} \label{gauto} Given $\bfe$ and $P$, there exists a finite automaton $\mathcal{F}$ over the alphabet $\mathcal{A}$ such that the exponent of $p$ in $\mathfrak{S}^\bfe_T(P)$ equals $\lambda(T)$ if and only if $\mathcal{F}$ accepts $T$. \end{lemma} 

The proof is essentially identical to that of Lemma \ref{auto}, and we omit it. We finish by lifting the non-negativity restriction which we imposed on $\bfe$. Below is the generalisation of Lemma \ref{neg}:

\begin{lemma} \label{gneg} Both Lemmas \ref{gdiv} and \ref{gauto} remain true when the components of $\bfe$ are allowed to be arbitrary (possibly negative) integers. \end{lemma} 

Once again, the proof proceeds along the same lines as that of the original Lemma \ref{neg}. We limit ourselves to highlighting a couple of steps where the generalisation is not entirely trivial. For all $i \ge 0$ and all integers $x$ not divisible by $p$, we know that $x^{(p - 1)p^i} \equiv 1 \pmod{p^{i + 1}}$ by virtue of $\varphi(p^{i + 1}) = (p - 1)p^i$ and Euler's theorem. Furthermore, in connection with the construction of the finite automaton $\mathcal{F}$, we must invoke part (a) of Lemma \ref{modp2}.

\section{Proof of the Theorem} \label{proof-i}

We now specialise the toolbox developed in Section \ref{abac-ii} to parts (c) and (d) of the Theorem. Our argument will be based on the sum ($\mathbf{W}^\mathsection$) of Section \ref{init}. We aim to show that the exponent of $2$ in ($\mathbf{W}^\mathsection$) is at least $2s_2(n) + \delta$; and we must also describe the cases of equality.

So we set $p = 2$ and $r = g = d/2$. Our alphabet will be $\mathcal{A} = \{\sfM_0, \sfM_1, \ldots, \sfM_g\}$, with the letter $\sfM_i$ corresponding to the multiset which consists of $i$ copies of the binary digit $1$ and $g - i$ copies of the binary digit~$0$. For convenience, we denote the number of occurrences of the letter $\sfM_i$ in the abacus type $T$ by $\alpha_i(T)$.

Let $\bfx = (x_1, x_2, \ldots, x_g)$ be any vector with nonnegative integer components which satisfies $n = x_1 + x_2 + \cdots + x_g$. Observe that, if $\bfx$ is associated with an abacus of type $T$, then $s_2(x_1) + s_2(x_2) + \cdots + s_2(x_g) = \sum_{i = 0}^g i\alpha_i(T)$ and $s_2(n) \le \sum_{i = 0}^g s_2(i)\alpha_i(T)$.

The first few steps of the argument will be the same as in the proof of Proposition \ref{dn}. Once again, we set $\bfe = (-4, 1)$ and $P = 1$. Then we split ($\mathbf{W}^\mathsection$) into sub-sums over abacus types. Fix an abacus type $T$, and consider the sub-sum which corresponds to it. Just as in the proof of part (b) of the Theorem, the exponent $\omega(\bfx)$ of $2$ in the $\bfx$-summand $E(\bfx) = \binom{n}{x_1, x_2, \ldots, x_g}^2 \binom{2x_1}{x_1} \binom{2x_2}{x_2} \cdots \binom{2x_g}{x_g}$ of ($\mathbf{W}^\mathsection$) equals $3(s_2(x_1) + s_2(x_2) + \cdots + s_2(x_g)) - 2s_2(n)$. Hence, $\omega$ is constant over $T$. Denote its value by $\omega(T)$. We get that $\omega(T) \ge \sum_{i = 0}^g (3i - 2s_2(i))\alpha_i(T)$.

We can now rewrite our sub-sum as \[\sum_T E(\bfx) = 2^{\omega(T)} \cdot \Theta_2(n!)^2 \cdot \mathfrak{S}^{(-4, 1)}_T(1).\]

We cannot invoke Lemma \ref{gneg} directly, as the presence of the constant term $\delta$ means that it is not strong enough anymore. Compared to the applications of our technique in the proofs of Propositions \ref{odd}--\ref{gdn}, this time around the argument will involve one extra subtlety.

Specifically, we are not going to halt our reductions when we reach the empty abacus type. Instead, we will halt them one step earlier than that, when we reach a single-letter abacus type. The point of this modification is as follows: Our desired exponent of $2$ is a linear function of~$s_2(n)$. Roughly speaking, the series of reductions we do carry out will account for the slope of this linear function. On the other hand, the residual sum over a single-letter abacus type which remains at the end will be responsible for the constant term.

Hence, we split $T$ into its first letter $\sfM_h$ and a ``tail'' $H$. We can assume without loss of generality that $h \neq 0$, because otherwise the leftmost column of every abacus of type $T$ will consist entirely of meaningless leading zeroes.

Fix a positive integer $k$ with $k \ge 2s_2(n) + \delta$. Then \[\mathfrak{S}^{(-4, 1)}_T(1) \equiv \mathfrak{S}^{(2^k - 4, 1)}_T(1) \pmod{2^{k + 2}}\] as in the proof of Lemma \ref{gneg}.

By the same reasoning as in the proof of Lemma \ref{gdiv}, we can find an integer vector $\bff$ and a symmetric integer-coefficient $g$-variable polynomial $R$ such that \[\mathfrak{S}^{(2^k - 4, 1)}_T(1) = 2^{\lambda(H)} \cdot \mathfrak{S}^\bff_{\sfM_h}(R).\]

Since it is symmetric, $R$ must be constant over $\varrho(\sfM_h)$. Denote its value by $R(\sfM_h)$. Then \[\mathfrak{S}^\bff_{\sfM_h}(R) = \binom{g}{h} \cdot R(\sfM_h) \cdot \Omega^\dagger\] for some odd positive integer $\Omega^\dagger$.

Putting the pieces together, we arrive at \[\sum_T E(\bfx) \equiv 2^{\omega(T) + \lambda(H)} \cdot \binom{g}{h} \cdot R(\sfM_h) \cdot \Omega^\ddagger \pmod{2^{k + 2}},\] with $\Omega^\ddagger = \Theta_2(n!)^2 \cdot \Omega^\dagger$ being an odd positive integer as well.

We will show next that the right-hand side is divisible by $2$ with exponent at least $2s_2(n) + \delta$. We are going to consider the cases of $h = 1$ and $h \ge 2$ separately.

Suppose first that $h \ge 2$. Then $\lambda(H) \ge \alpha_1(T)$. By our earlier estimates for $s_2(n)$ and $\omega(T)$ in terms of $\alpha_0(T)$, $\alpha_1(T)$, $\ldots$, $\alpha_g(T)$, it follows that \[\omega(T) + \lambda(H) \ge 2s_2(n) + \sum_{i = 2}^g (3i - 4s_2(i)) \cdot \alpha_i(T).\]

However, it is routine to check that $3x > 4s_2(x)$ for all positive integers $x \ge 2$. Hence, $\omega(T) + \lambda(H) \ge 2s_2(n)$; and the slope of our linear function of $s_2(n)$ has been accounted for.

We move on to the constant term. Since $2 \le h \le g$ and $\alpha_h(T) \ge 1$, we may ``borrow'' $3h - 4s_2(h)$ out of the preceding estimate. Furthermore, by Kummer's theorem, the exponent of $2$ in $\binom{g}{h}$ equals $s_2(h) + s_2(g - h) - s_2(g)$. We claim that the sum of the latter couple of quantities is in fact strictly greater than $\delta$.

Indeed, $\delta = \nu_2(g) - 1$ and $\nu_2(g) = s_2(g - 1) + 1 - s_2(g)$. (One way to verify this identity is to take a look at the binary expansions of $g - 1$ and $g$; another is to apply Kummer's theorem to the binomial coefficient $\binom{g}{1}$.) So we can rewrite our desired inequality as $s_2(g - h) + 3h > s_2(g - 1) + 3s_2(h)$. The case of $h = 2$ is clear. Otherwise, if $h \ge 3$, we set $x = h - 1$ in our earlier observation that $3x > 4s_2(x)$ whenever $x \ge 2$. We add together the resulting inequality and $s_2(g - h) + s_2(h - 1) \ge s_2(g - 1)$ so as to obtain that $s_2(g - h) + 3h > s_2(g - 1) + 3s_2(h - 1) + 3$. The final step follows by the obvious $s_2(h - 1) + 1 \ge s_2(h)$.

This confirms our claim; and so the constant term of our linear function of $s_2(n)$ has been settled as well. We conclude that, if $h \ge 2$, then the exponent of $2$ in the sub-sum $\sum_T E(\bfx)$ is in reality strictly greater than what we require.

We continue with the case when $h = 1$. Then $\lambda(H) \ge \alpha_1(T) - 1$ and, reasoning as above, we find that $\omega(T) + \lambda(H) \ge 2s_2(n) - 1$. Furthermore, the exponent of $2$ in $\binom{g}{1}$ becomes $\nu_2(g) = \delta + 1$. Thus in this case, too, things work out as they should.

What remains is to sort out the equality conditions. We can tell right away that the exponent of $2$ in the sub-sum $\sum_T E(\bfx)$ will be strictly greater than $2s_2(n) + \delta$ unless $h = 1$ and $\alpha_i(T)$ vanishes for all $2 \le i \le g$. Or, in other words, we may safely focus on the abacus type $T = \cf(n)$. Our analysis shows also that the exponent of $2$ in the corresponding sub-sum $\sum_{\cf(n)} E(\bfx)$ will be exactly equal to $2s_2(n) + \delta$ if and only if $R(\sfM_1) = R(0, 0, \ldots, 0, 1)$ is odd. To determine the parity of this number, we construct a suitable finite automaton as in the proof of Lemma \ref{gneg}.

\begin{table}[ht] \renewcommand{\arraystretch}{1.125} \footnotesize \null \hfill \parbox{0.475\textwidth}{\centering \begin{tabular}{|c|c|c|}
from & with & to\\
\hline
\multirow{2}{*}{$1$} & $0$ & $1$\\
\cline{2-3}
& $1$ & $\sum x_i + 1$\\
\hline
\multirow{2}{*}{$\sum x_i + 1$} & $0$ & $1$\\
\cline{2-3}
& $1$ & $0$\\
\hline
$0$ & $0$, $1$ & $0$\\
\hline
\end{tabular} \caption{} \label{4mod8}} \hfill \parbox{0.475\textwidth}{\centering \begin{tabular}{|c|c|c|}
from & with & to\\
\hline
\multirow{2}{*}{$1$} & $0$ & $1$\\
\cline{2-3}
& $1$ & $\sum x_i$\\
\hline
\multirow{2}{*}{$\sum x_i$} & $0$ & $0$\\
\cline{2-3}
& $1$ & $\sum x_i$\\
\hline
$0$ & $0$, $1$ & $0$\\
\hline
\end{tabular} \caption{} \label{0mod8}} \hfill \null \end{table}

The calculations are straightforward. We use the alphabet $\{0, 1\}$ in place of $\mathcal{A}$ because the binary expansion of $n$ maps onto $\cf(n)$ by means of the substitutions $0 \to \sfM_0$ and $1 \to \sfM_1$. The ``standard'' construction prescribes that each state must be labelled with an ordered pair consisting of a fold and a remainder. However, the folds turn out not to influence the behaviour of the states, and so we omit them. When $\nu_2(d) = 2$, we obtain the finite automaton of Table~\ref{4mod8}. Otherwise, when $\nu_2(d) \ge 3$, we obtain the finite automaton of Table~\ref{0mod8} instead. The desired description of all $n$ which attain equality is immediate, and our proof of parts (c) and (d) of the Theorem is complete.

\section{Alternative Proof of the Theorem} \label{proof-ii}

The identity \[W^\star_4(n) = \sum_{x + y = n} (-1)^x\binom{2x}{x}\binom{3x + y}{x, x, x, y}2^{4y} \tag{\textbf{X}}\] for the $n$-th Domb number, of Chan and Zudilin \cite{CZ}, makes it quite easy to analyse $w^\star_4(n)$. Indeed, almost all summands on the right-hand side are divisible by very high powers of $2$, and so it is only a tiny number of summands that actually matter for the analysis. Here, we present an alternative proof of parts (c) and (d) of the Theorem rooted at this observation.

The full ``derivation tree'' will be as follows: First, we use (\textbf{X}) to resolve $d = 4$. (Lemma \ref{alt-1}.) Then, out of this special case, we get all $\nu_2(d) = 2$. (Lemma \ref{alt-2}.) Once again out of $d = 4$, by means of a different argument, we get $d = 8$ as well. (Lemma \ref{alt-3}.) Finally, out of $d = 8$, we derive all $\nu_2(d) \ge 3$. (Lemma \ref{alt-4}.)

It is instructive to compare the two proofs. The new one, on the face of it, would seem simpler. However, this simplicity is deceptive; a lot of the difficulty is hidden away in the ``magic'' identity~(\textbf{X}). Even if we set aside the question of its verification, we must still ask: How does one come up with (\textbf{X}), in a motivated manner, while thinking about high-dimensional lattice walks?

One more consideration is that the new proof feels somewhat ``ad hoc'', with different tricks employed for different values of $d$. By contrast, our previous proof resolves all values of $d$ (in parts (c) and (d) of the Theorem) uniformly, by means of the same general method.

We return now to the task at hand. Our quick overview of the logical structure of the argument shows that we will often be reducing bigger values of $d$ to smaller ones. This is accomplished with the help of the identity ($\mathbf{W}^\square$) given below. Suppose that $d = ab$. Then \[W_d(n) = \sum_{x_1 + x_2 + \cdots + x_b = n} \binom{2n}{2x_1, 2x_2, \ldots, 2x_b} W_a(x_1)W_a(x_2) \cdots W_a(x_b). \tag{$\mathbf{W}^\square$}\]

Indeed, recall the definition of the Laurent polynomial $S_d$ from the introduction. Let $\bfz = (z_1, z_2, \ldots, z_d)$. For each $i = 1$, $2$, $\ldots$, $b$, let also $\bfz_i = (z_{a(i - 1) + 1}, z_{a(i - 1) + 2}, \ldots, z_{ai})$. Then $S_d(\bfz) = S_a(\bfz_1) + S_a(\bfz_2) + \cdots + S_a(\bfz_b)$. Raising both sides to the power of $2n$, we find that \[S_d(\bfz)^{2n} = \sum_{e_1 + e_2 + \cdots + e_b = 2n} \binom{2n}{e_1, e_2, \ldots, e_b}S_a(\bfz_1)^{e_1}S_a(\bfz_2)^{e_2} \cdots S_a(\bfz_b)^{e_b}.\]

Since the polynomials $S_a(\bfz_1)$, $S_a(\bfz_2)$, $\ldots$, $S_a(\bfz_b)$ do not share any variables, the constant term of $S_a(\bfz_1)^{e_1}S_a(\bfz_2)^{e_2} \cdots S_a(\bfz_b)^{e_b}$ equals the product of the constant terms of $S_a(\bfz_1)^{e_1}$, $S_a(\bfz_2)^{e_2}$, $\ldots$, $S_a(\bfz_b)^{e_b}$. Since, furthermore, the constant term of $S_a(\bfz_i)^{e_i}$ is nonzero if and only if $e_i$ is even, the desired conclusion now follows. We are all set up to begin with the proof.

\begin{lemma} \label{alt-1} The Theorem holds when $d = 4$. \end{lemma}

\begin{myproof} We must show that the exponent of $2$ in the right-hand side of (\textbf{X}) is always at least $2s_2(n)$, and also describe the equality cases.

By Kummer's theorem, the exponent of $2$ in $\binom{2x}{x}$ is $s_2(x)$ and in $\binom{3x + y}{x, x, x, y}$ it is $3s_2(x) + s_2(y) - s_2(3x + y)$. Since $2s_2(x) + s_2(y) = s_2(3)s_2(x) + s_2(y) \ge s_2(3x + y)$, the latter quantity cannot fall below $s_2(x)$. Of course, also $y \ge s_2(y)$. Hence, the exponent of $2$ in the $(x, y)$-summand becomes at least $2s_2(x) + 4s_2(y) = 2(s_2(x) + s_2(y)) + 2s_2(y) \ge 2s_2(n) + 2s_2(y)$. This is strictly greater than $2s_2(n)$ unless $x = n$ and $y = 0$.

We are left to analyse the $(n, 0)$-summand. The chain of inequalities in the preceding paragraph, specialised with $x = n$ and $y = 0$, shows that the exponent of $2$ in it will be exactly equal to $2s_2(n)$ if and only if the multiplication of $3$ and $n$ is carry-free in binary. Clearly, this condition is equivalent to the binary expansion of $n$ not containing two adjacent~$1$s. \end{myproof}

\begin{lemma} \label{alt-2} The Theorem holds when $\nu_2(d) = 2$. \end{lemma}

\begin{myproof} We apply ($\mathbf{W}^\square$) with $a = 4$ and $b = d/4$. Notice that we may discard all summands on the right-hand side of ($\mathbf{W}^\square$) which are divisible by $2$ with exponent at least $3s_2(n) + 1$.

By Lemma \ref{alt-1}, the exponent of $2$ in the product $W_4(x_1)W_4(x_2) \ldots W_4(x_b)$ is at least $3(s_2(x_1) + s_2(x_2) + \cdots + s_2(x_b))$. If $x_1$, $x_2$, $\ldots$, $x_b$ do not form a carry-free partitioning of $n$, then $s_2(x_1) + s_2(x_2) + \cdots + s_2(x_b) \ge s_2(n) + 1$ and the associated summand of ($\mathbf{W}^\square$) can be safely discarded.

Suppose, from now on, that $x_1$, $x_2$, $\ldots$, $x_b$ do form a carry-free partitioning of $n$. Then the multinomial coefficient $\binom{2n}{2x_1, 2x_2, \ldots, 2x_b}$ in the associated summand of ($\mathbf{W}^\square$) becomes odd. What remains is to describe those $n$ which admit an odd number of carry-free partitionings where each piece satisfies the equality conditions of Lemma \ref{alt-1}.

Suppose that, in the binary expansion of $n$, there are $k$ blocks of $1$s of lengths $c_1$, $c_2$, $\ldots$, $c_k$, respectively. (Here, by a ``block'' of $1$s we mean a maximal continuous run of $1$s.) For each block of length $c$, the number of ways to distribute its $1$s between $x_1$, $x_2$, $\ldots$, $x_b$ so that no $x_i$ gets two adjacent $1$s works out to $b(b - 1)^{c - 1}$. (Say we assign the $1$s of our block to $x_1$, $x_2$, $\ldots$, $x_b$ one by one, going from left to right. There will be $b$ options to choose from on the first step and $b - 1$ options to choose from on each one of the subsequent $c - 1$ steps.) Since $b = d/4$ is odd, the product of these numbers over all $c = c_1$, $c_2$, $\ldots$, $c_k$ will be odd as well if and only if $c_1 = c_2 = \cdots = c_k = 1$, as desired. \end{myproof}

\begin{lemma} \label{alt-3} The Theorem holds when $d = 8$. \end{lemma}

\begin{myproof} We apply ($\mathbf{W}^\square$) with $a = 4$ and $b = 2$. Notice that we may discard all summands on the right-hand side of ($\mathbf{W}^\square$) which are divisible by $2$ with exponent at least $3s_2(n) + 2$. Just as in the proof of Lemma \ref{alt-2}, this includes all summands where $x_1$ and $x_2$ do not form a carry-free partitioning of $n$.

Observe next that $x_1 \neq x_2$ in all carry-free partitionings of $n$. So, by swapping $x_1$ and $x_2$, we can pair up all summands on the right-hand side of ($\mathbf{W}^\square$) which correspond to such partitionings. Each pair will consist of two copies of the same positive integer. Furthermore, in each pair, the corresponding binomial coefficient $\binom{2n}{2x_1, 2x_2}$ will be odd. What remains is to determine when the number of carry-free partitionings of $n$ such that both of $x_1$ and $x_2$ satisfy the equality conditions of Lemma \ref{alt-1} is twice an odd positive integer.

Define $c_1$, $c_2$, $\ldots$, $c_k$ as in the proof of Lemma \ref{alt-2}. By the same reasoning as before, a block of length $c$ can be distributed between $x_1$ and $x_2$ in $2 \cdot 1^{c - 1} = 2$ ways. (The distribution must alternate between $x_1$ and $x_2$ in a zig-zag pattern.) Hence, the total number of partitionings which satisfy our constraints works out to $2^k$. This is twice an odd positive integer if and only if $k = 1$, as desired. \end{myproof}

\begin{lemma} \label{alt-4} The Theorem holds when $\nu_2(d) \ge 3$. \end{lemma}

\begin{myproof} We apply ($\mathbf{W}^\square$) with $a = 8$ and $b = d/8$ -- but there is one new subtlety this time around. We are permitting some of $x_1$, $x_2$, $\ldots$, $x_b$ to vanish; while the Theorem and Lemmas \ref{alt-1}--\ref{alt-4} all assume that $n$ is a positive integer. This did not matter much in the proofs of Lemmas \ref{alt-2} and \ref{alt-3} because Lemma \ref{alt-1} remains true when $n = 0$. However, Lemma \ref{alt-3}, which we are about to use now, does not.

Clearly, permuting $x_1$, $x_2$, $\ldots$, $x_b$ does not alter the value of the corresponding summand. We sort all summands on the right-hand side of ($\mathbf{W}^\square$) into classes so that two summands belong to the same class if and only if they can be obtained from one another by means of such permutations. For each multiset $X = \{x_1, x_2, \ldots, x_b\}$, let $f(X)$ be the number of summands in the class of $X$ and let $F(X)$ be the value of each such summand.

Denote by $h$ the number of non-zeroes in $X$. Then $f(X)$ is a multiple of $\binom{b}{h}$ because this is the number of ways to choose which pieces of the partitioning are going to be assigned a nonzero element of $X$; and the exponent of $2$ in this binomial coefficient equals $s_2(h) + s_2(b - h) - s_2(b)$ by Kummer's theorem. Furthermore, $F(X)$ is divisible by $2$ with exponent at least $3(s_2(x_1) + s_2(x_2) + \cdots + s_2(x_b)) + h$ by Lemma \ref{alt-3}. We conclude that the total contribution $f(X)F(X)$ of $X$ towards the right-hand side of ($\mathbf{W}^\square$) is divisible by $2$ with exponent which cannot fall below the sum of the latter couple of quantities.

Once again, we may discard all classes where this sum is strictly greater than $3s_2(n) + \delta$. Of course, $3(s_2(x_1) + s_2(x_2) + \cdots + s_2(x_b)) \ge 3s_2(n)$. We claim that, if $h \ge 3$, additionally $s_2(h) + s_2(b - h) - s_2(b) + h > \delta$.

Indeed, $\delta = \nu_2(b) + 1$ and $\nu_2(b) = s_2(b - 1) + 1 - s_2(b)$. (We have already encountered this particular way of expressing $\nu_2$ in terms of $s_2$ in Section \ref{proof-i}.) Since also $s_2(b - h) + s_2(h - 1) \ge s_2(b - 1)$, the problem boils down to verifying the inequality $s_2(h) - s_2(h - 1) + h \ge 3$. Or, equivalently, $h - \nu_2(h) \ge 2$; which is obviously true when $h \ge 3$. This confirms our claim.

Thus we can safely discard all classes where $h \ge 3$. Suppose, next, that $h = 2$. Then our estimates in the preceding couple of paragraphs yield a lower bound of $3s_2(n) + \delta$ for the exponent of $2$ in $f(X)F(X)$. However, if the two nonzero elements of $X$ are equal, they cannot form a carry-free partitioning of $n$; whereas, if they are distinct, $f(X) = 2\binom{b}{2}$ is divisible by $2$ with a strictly greater exponent than $\binom{b}{2}$. Either way, our lower bound of $3s_2(n) + \delta$ cannot be attained.

What remains is to look into the unique class where $h = 1$ and $X = \{0, 0, \ldots, 0, n\}$. Its contribution towards the right-hand side of ($\mathbf{W}^\square$) amounts to $bW_8(n)$. Since $\nu_2(b) = \nu_2(d/8) = \delta - 1$, the desired result now follows immediately by Lemma \ref{alt-3}. \end{myproof}

This completes our alternative proof of parts (c) and (d) of the Theorem.

\section{Further Work} \label{further}

Our initial hunch that the exponents of $2$ in the numbers $W_d(n)$ might be worth studying was based on a connection to the numbers $B(m, n)$. We go on now to discuss this connection in somewhat deeper detail.

Consider a grid of size $2m \times 2n$ whose cells are coloured in black and white. The colouring is \emph{balanced} when, in each row and each column, half of the cells are white and half are black. We denote the number of such colourings by $B(m, n)$. The conjecture of Bhattacharya cited in the introduction states that the exponent of $2$ in $B(m, n)$ equals $s_2(m)s_2(n)$. A few partial results are collected in \cite{B}. However, a full proof has so far been out of reach.

We may view $B(m, n)$ and $W_d(n)$ as two special cases of a single overarching counting problem. Given two nonnegative integers $k$ and $\ell$, we write $U_{k, \ell}(n)$ for the number of ways to colour the cells of a $(k + \ell) \times 2n$ grid in black and white so that each row (of size $k + \ell$) contains $k$ cells of one colour and $\ell$ of the other, while each column (of size $2n$) is balanced. Then $B(m, n) = U_{m, m}(n)$ and, if $d \neq 2$, also $W_d(n) = U_{1, d - 1}(n)$. (When $d = 2$, the relation becomes $W_2(n) = \binom{2n}{n}U_{1, 1}(n)$ instead.)

Notice that $U_{k - 1, k}(n) = U_{k, k}(n)$. We can now recover the relations $W_1(n) = U_{1, 0}(n) = U_{1, 1}(n) = B(1, n)$ and $W_3(n) = U_{1, 2}(n) = U_{2, 2}(n) = B(2, n)$ referenced in the introduction.

Does our analysis of high-dimensional lattice walks throw any light on balanced grid colourings? For the author, a crucial part of the appeal of the former problem was the likelihood that any methods developed for it might turn out to be helpful with the latter one, too. We round off the paper by sketching, in rather broad strokes, some ideas in this vein.

Our proof of Lemma \ref{alt-1} in Section \ref{proof-ii} relies on a representation of the desired count in the form $\sum_{i = 0}^n 2^{ei}F(i)$ with an ``arithmetically well-behaved'' $F$. What are the positive integers $m$ such that $B(m, n)$, considered as a function of $n$, admits a similar representation? Of course, we are looking specifically for representations which would allow us to determine the exponent of $2$ in $B(m, n)$ by examining just a small number of summands near the beginning.

On the other hand, our proof of Proposition \ref{dn} in Section \ref{app} (whose statement is, in fact, equivalent to that of the aforementioned Lemma \ref{alt-1}) uses a different representation of the same count in which a lot of the summands are divisible by $2$ with exponents much smaller than the one we seek. However, by bundling these summands together into sub-sums over abacus types, we were able to boost the divisibility as needed. Do the domino representations of $B(m, n)$ constructed in Sections 4 and 5 of \cite{B} admit similar boosts?

One more route to explore would be as follows: We have already learned quite a bit about the exponents of $2$ in the numbers $U_{1, \ell}(n)$. We could move on next to the numbers $U_{2, \ell}(n)$. Then $U_{3, \ell}(n)$; and so on, and so forth. We could thus build up our understanding of the numbers $U_{k, \ell}(n)$ little by little, via a series of problems which gradually increase in difficulty. Any such general understanding would naturally illuminate the numbers $B(m, n)$ as well.

\section*{Acknowledgements}

The present paper was written in the course of the author's PhD studies under the supervision of Professor Imre Leader. The author is thankful to Prof.\ Leader for his unwavering support.


\begin{thebibliography}{99}

\bibitem{O1} OEIS entry A287318.

\bibitem{O2} OEIS entry A287316.

\bibitem{O3} OEIS entry A002895.

\bibitem{C} Neil Calkin, Factors of Sums of Powers of Binomial Coefficients, \emph{Acta Arithmetica}, 1998.

\bibitem{RS} L.\ Bruce Richmond and Jeffrey Shallit, Counting Abelian Squares, \emph{The Electronic Journal of Combinatorics}, 2009.

\bibitem{CZ} Heng Huat Chan and Wadim Zudilin, New Representations for Ap\'{e}ry-Like Sequences, \emph{Mathematika}, 2010.

\bibitem{OS} Robert Osburn and Brundaban Sahu, A Supercongruence for Generalized Domb Numbers, \emph{Functiones et Approximatio Commentarii Mathematici}, 2013.

\bibitem{D} \'{E}ric Delaygue, Arithmetic Properties of Ap\'{e}ry-Like Numbers, arXiv preprint, 2013; \emph{Compositio Mathematica}, 2018.

\bibitem{P} Bjorn Poonen, personal communication, 2025.

\bibitem{B} Nikolai Beluhov, Powers of 2 in Balanced Grid Colourings, \emph{Enumerative Combinatorics and Applications}, 2026.

\end{thebibliography}
\end{document}